\documentclass[11pt]{article}
\usepackage{amsfonts}
\usepackage{latexsym,amsmath}
\usepackage{array}
\usepackage{amssymb}  % for \square
\usepackage[symbol]{footmisc}
\usepackage{fancyhdr}
\usepackage[a4paper, portrait, margin=1.1811in]{geometry}
\usepackage{amsthm}
\usepackage{relsize}
\usepackage{graphicx}
\usepackage{lscape}
\usepackage{booktabs}
\usepackage{longtable}
\usepackage{multirow}
\usepackage{threeparttable}
\usepackage[utf8]{inputenc}
\usepackage[T1]{fontenc}
\usepackage[table]{xcolor}
\usepackage{booktabs}
\usepackage{tabularx}
\arrayrulecolor{black}
\usepackage[longnamesfirst,sort]{natbib}% Citation support using natbib.sty
\bibpunct[, ]{(}{)}{;}{a}{,}{,}% Citation support using natbib.sty

\usepackage[colorlinks, citecolor=cyan]{hyperref}
\newtheorem{t1}{Theorem}[section]

\newtheorem{d1}{Definition}[section]

\newtheorem{ex}{Example}[section]

\renewcommand{\arraystretch}{1.3}
\begin{document}
	\title{\textbf{Quantile-based Fractional Generalized Cumulative Past Entropy}}
	\author{Poulami Paul\and Chanchal Kundu\footnote{\textit{Corresponding
				author e-mail}: \color{cyan}ckundu@rgipt.ac.in;
			\color{cyan}chanchal$_{-}$kundu@yahoo.com.}
		\and Department of Mathematical Sciences\\
		Rajiv Gandhi Institute of Petroleum Technology\\
		Jais 229 304, U.P., India}
	\date{September, 2025}
	\maketitle
	\begin{abstract}	
		Uncertainty in past lifetime distributions and the timing of inactivity in systems and their components can be effectively measured using the fractional generalized cumulative past entropy (FGCPE) and its dynamic extension (DFGCPE), introduced by Di Crescenzo et al. (2021). Building on this framework, we propose a quantile-based variant, the quantile fractional generalized cumulative past entropy (QFGCPE), along with its dynamic time-dependent counterpart (DQFGCPE).
		 Closed-form expressions of these measures are derived analytically for a variety of lifetime distributions, including those with and without explicit distribution functions. Fundamental properties such as bounds, monotonicity, and stochastic orderings are investigated to assess robustness and interpretability. Furthermore, we construct a nonparametric estimator of QFGCPE and establish its asymptotic validity through extensive simulation studies involving bias, mean squared error (MSE), and root mean squared error (RMSE). Finally, the sensitivity of the proposed QFGCPE measure is examined by comparing its behavior with the logistic map, demonstrating its ability to capture transitions from order to chaos.

	\end{abstract}	
	{\bf Key Words and Phrases:} Fractional generalized cumulative entropy, Quantile form, Stochastic ordering, Nonparametric estimation, Logistic map \\
	{\bf MSC2020 Classifications:} Primary 94A17; Secondary 62B10, 37M25, 62G05, 60E15. %, 65C05.
	
	\section{Introduction} 
	The uncertainty in a probabilistic system is a key determinant of its overall efficiency. Assessing the system uncertainty is fundamental to risk and survival analysis as well as to problems pertaining to statistical planning and inference. Typically, while quantifying the uncertainty of a system, the distribution function of the system components or random variables (r.vs) such as the returns of a financial portfolio asset are to be determined. However, this process of characterizing every random variables by a distribution function (DF) is not that simple, since there exist certain r.vs that do not admit tractable distribution functions. Examples include the Govindarajalu distribution \citep{go}, various forms of lambda distributions \citep{vl} and the power-Pareto distribution \citep{hl}. In such cases, quantile functions (QFs) provide a practical alternative for describing the system and predicting its uncertainty by working with the inverse behavior of the cumulative distribution.
	
	At its core, the entropy initially designed by \cite{s}, being a function of DF,
	 serves as a fundamental measure of uncertainty and information within data. For a non-negative absolutely continuous r.v $X$ with DF $G(x)$ and probability density function (pdf) $g(x)$, the Shannon differential entropy is defined as:
	\begin{equation}
	   H_S(X) = -\int_{0}^{\infty}g(x) \ln g(x), \qquad g(x) > 0. 
	   \label{Shannon}
	\end{equation} This entropy measure \eqref{Shannon} plays a central role in information theory, reliability analysis, and statistical modeling. 
	Since the pioneering work of Shannon, various extensions and generalizations of entropy have been proposed to capture different aspects of randomness and information content in complex systems. 
	Among these measures, cumulative entropies have gained significant attention due to their simplicity and the ease with which they capture distributional information via empirical cumulative and survival functions, thereby avoiding the complexity of density function estimation and yielding tractable yet refined measures particularly suited to reliability and lifetime analysis \citep[cf.][and the references therein]{dl,r,dt,nv}.
	The cumulative entropy of $X$ with survival function (SF) $ \overline{G}(x) $, as proposed by \cite{r}, is defined as:
	\begin{equation}
		CR\xi(X) = -\int_{0}^{\infty}\overline{G}(x) \ln \overline{G}(x)
		\label{1.2}
	\end{equation}
	and the cumulative past entropy of $X$ with DF $G(x)$, as formulated by \cite{dl}, is expressed as:
	\begin{equation}
		CP\xi(X) = -\int_{0}^{\infty}G(x) \ln G(x).
		\label{1.3}
	\end{equation}
	 % ; Navarro et al., 2010). %\citep{DiCrescenzo2009, Navarro2011}.
	
	These entropy measures are further generalized through the use of QF in place of DF in its definitions since QFs are independent of the distribution of random variables \citep{ss,ssj}. Moreover, the QFs give equal importance to all the sample points, including the outliers, unlike DFs, thereby making them potentially more able in giving good reliability analysis results even with small sample size or limited information \citep[cf.][]{nsb,ns}. Several other properties of QFs are studied which proves it to be an useful tool for statistical analysis and model identification \citep[cf.][and the references therein]{gi,ne,nva,vns,vh,as,kba} %(cf. Gilchrisht, 2000; Nair et al., 2013, Vineshkumar and Nair, 2021, Varkey and Haridas, 2023, Aswin et al., 2023, Kayal and Balakrishnan, 2024 ).
	%The quantile based cumulative entropies are given as:
	%The quantile-based definitions can help in describing r.vs even with intractable DFs and therefore, are becoming more popular among researchers in dealing with non-parametric data. 
	
	While traditional methods often focus on mean-based predictions, they frequently fall short in capturing the full spectrum of distributional characteristics, especially in volatile or complex systems. This is where the power of QF-based methods becomes indispensable. Quantile-based approaches allow us to model and predict various points across the conditional distribution, not just the average. This provides a richer understanding of data behavior, enabling us to quantify risks and
	uncertainties more precisely. When combined with entropy, these techniques offer unparalleled insights into the probabilistic nature of physical phenomena, revealing hidden patterns and potential extremes that might otherwise be overlooked.
	
	Let us consider a non-negative absolutely continuous random variable $ X $ identifying the lifetime of a component (or individual) or the whole system (or population) with cumulative distribution function (cdf) $ G $ and probability density function (pdf) $ g. $ Then the QF can be defined as the inverse of the cdf given as:
	\begin{equation}
		Q(v) = G^{-1}(v) = \inf\{x|G(x)\geq v\}, 0\leq v \leq 1.
		\label{q1}
	\end{equation} 
	Here, $ g(Q(v)) $ or $gQ(v)$ and $ q(v)= Q^{'}(v) = \frac{dQ(v)}{dv} $ are known as the density quantile function (dqf) and the quantile density function (qdf), respectively \citep[cf.][]{p}.
	Hence, on differentiating both sides of $ G(Q(v)) ~\text{or}~ GQ(v) = v $ obtained from (\ref{q1}) with respect to $ v $, we have that 
	\begin{equation}
		q(v)g(Q(v)) = 1.
		\label{q2}
	\end{equation}
	Then the quantile-based cumulative residual $CR\xi_Q$ and cumulative past entropy $CP\xi_Q$ can be derived from Eqs. \eqref{1.2} and \eqref{1.3}, respectively, using Eqs. \eqref{q1} and \eqref{q2} as \cite{ssj}: 
	\begin{equation}
		CR\xi_Q(X) = -\int_{0}^{1}(1-p) \ln (1-p) q(p) dp
	\end{equation}
	and 
	\begin{equation}
		CP\xi_Q(X) = -\int_{0}^{1}p \ln (p) q(p) dp,
	\end{equation} respectively.
	
	%\hspace*{0.2in}
	%Quantile-based cumulative entropy measures provide a framework for analyzing uncertainty in random variables using quantile functions, particularly useful when distribution functions are intractable.   
	
	 Recently, \cite{xi} and \cite{dc} %Di Crescenzo et al.\ (2021)
	 introduced the concepts of fractional cumulative residual entropy (FCRE) and fractional generalized cumulative past entropy (FGCPE), respectively, by replacing the density function in the fractional entropy definition proposed by \cite{u} with SF and DF, respectively. These entropy measures were formulated to find a flexible extension that combines fractional calculus with the cumulative entropies defined earlier by \cite{r} and \cite{dl}, respectively, to better capture structural properties of lifetimes and reliability functions. This framework has motivated several works exploring properties, bounds, and applications of fractional cumulative entropies in stochastic modeling and survival analysis \citep{kb,ks}. 
	 The FCRE of $X$ is defined as \citep{xi}:
	 \begin{equation}
	 	\mathcal{C}R\xi^\eta (X) = \int_{0}^{\infty} \overline{G}(x)[-\ln \overline{G}(x)]^\eta dx, ~ 0 \leq \eta \leq 1,
	 	\label{fcre}
	 \end{equation} and
	the FGCPE of $ X $ is defined as \citep{dc}:  %with a general non-negative weight function $ k(x) (\geq 0) $ 
	\begin{equation}
		\mathcal{C}P\xi^\eta (X) = \frac{1}{\Gamma (\eta +1)}\int_{0}^{\infty} G(x)[-\ln G(x)]^\eta dx, ~ \eta > 0,
		\label{fgcpe}
	\end{equation}
	provided the integral is finite, where $ \Gamma(\cdot) $ represents the gamma function.
	
	Motivated by the established advantages of the QFs $Q(\cdot)$ in (\ref{q1}), particularly their robustness and effectiveness in statistical inference relative to DFs, the quantile alternative of FCRE (QFCRE) of $X$ given in \eqref{fcre} is expressed as \citep{ses}:
	\begin{equation}
		\mathcal{C}R\xi_Q^\eta (X) = \int_{0}^{1} (1-p)[-\ln (1-p)]^\eta q(p) dp, \qquad 0 \leq \eta \leq 1,
		\label{qfcre}
	\end{equation} This quantile-based definition of FCRE, as given in Eq.~\eqref{qfcre}, facilitates the study of its statistical properties and enables a more straightforward construction of non-parametric estimators, as well as more efficient simulated and real data analyses. Furthermore, the quantile representation \eqref{qfcre} not only provided an alternative perspective on distributional uncertainty but also facilitated tractable closed-form analysis for a broad class of lifetime models.
	
	 % has been shown to be a better alternative to the DF in statistical modeling and data analysis by many authors, the quantile alternative of FCRE (QFCRE) has been proposed by \cite{ses}. 
	 %However, existing approaches are primarily density- or distribution-based, and quantile-based analogues are underdeveloped despite their recognized  (Parzen, 1979; Koenker2005). %\citep{Parzen1979, Koenker2005}.  
	
	In this paper, we propose and investigate the \emph{quantile-based fractional generalized cumulative past entropy} (QFGCPE), a novel quantile analogue of the FGCPE expressed by \eqref{fgcpe}. The remainder of this paper is organized into the following sections. Section 2 introduces the definition of the quantile version of FGCPE and examines some of its fundamental properties and bounds. Some important stochastic ordering relations are established and explicit closed-form expressions are obtained for several important distributions, including the uniform, exponential, Fréchet, half-logistic, power, Govindarajalu, and Davies distributions.
	In Section 3, we further extend the framework by introducing the dynamic version of Q-FGCPE, thereby enabling time-dependent uncertainty assessment in evolving systems. In Section 4, we develop a nonparametric estimator of Q-FGCPE from a methodological standpoint and study its finite-sample performance through Monte Carlo simulation. The empirical results, evaluated in terms of bias, mean squared error (MSE) and root mean squared error (RMSE) across varying sample sizes, demonstrate consistency and validate the estimator’s asymptotic properties. In Section 5, we further compare its behavior with the logistic map, a chaos indicator for classical dynamical systems	to strengthen the practical validity of the proposed measure, thereby highlighting the sensitivity of QFGCPE to transitions between regular and chaotic regimes. Finally, Section 6 concludes the study by portraying the important highlights and contributions of this paper.
	
	\section{Quantile-based FGCPE}
	Following the relation between the quantile function $Q$ and a continuous distribution function $G$ given by Eqs. \eqref{q1} and \eqref{q2}, we can define the quantile analogue of FGCPE as follows:
	\begin{d1}
		If we consider a non-negative continuous random variable $ X $ with quantile function $ Q(v), $ then the quantile-based FGCPE (or QFGCPE) is defined as:
		\begin{eqnarray}
			\mathcal{C}P\xi_Q^\eta &=& \frac{1}{\Gamma(\eta+1)}\int_{0}^{1} p(-\ln p)^\eta dQ(v) \nonumber\\
			&=& \frac{1}{\Gamma(\eta+1)}\int_{0}^{1} p(-\ln p)^\eta q(p) dp,  \eta > 0, \label{2.i}
		\end{eqnarray}
		where $ q(v) = \frac{d}{dv}Q(v) $ is the quantile density function (qdf).
	\end{d1}
	For a continuous random variable $ X $ with cdf $ G(x) $, pdf $ g(x), $ the reversed hazard rate function is:
	\begin{equation*}
		r(x) = \frac{g(x)}{G(x)}
	\end{equation*}
	Therefore the reversed hazard quantile function (RHQF) can be defined as \citep{nf,sj,q}:
	\begin{equation}
		R(v) = r(Q(v)) = rQ(v) = \frac{gQ(v)}{v} = [vq(v)]^{-1}.
		\label{RHQ}
	\end{equation}
	Hence, the QFGCPE represented by (\ref{2.i})can be expressed in terms of RHQF as
	\begin{equation}
		\mathcal{C}P\xi_Q^\eta = \int_{0}^{1} [R(p)]^{-1}(-\ln p)^\eta dp.
	\end{equation}
	The quantile functions and the respective QFGCPE for some important families of distributions are provided in Table \ref{tab1}.
	\begin{table}[!ht]
		\centering
		\caption{Quantile functions and $ \mathcal{C}P\xi_Q^\eta $ of some important lifetime distributions.}
		%\begin{tabular}{|l l l l l|}
		\renewcommand{\arraystretch}{1.3}
		\begin{tabular}{l l l l l }
			\toprule
			\textbf{Distribution} & \textbf{Parameters}  &  $ \mathbf{q(v)} $  &  $ \mathbf{Q(v)} $  &  $ \mathbf{\mathcal{C}P\xi_Q^\eta} $  \\
			\midrule
			Uniform  &  $[0,b], b > 0 $  &  $ b $  &  $ bv $  & $\frac{1}{2^{\eta+1}}$ \\
			Exponential  &  $\lambda > 0 $  & $ \frac{1}{\lambda(1-v)} $  &  -$\frac{1}{\lambda} \ln(1-v)$  &  $\frac{1}{\lambda} \left[\zeta(\eta+1)-1\right]$ \\
			Power  &  $0<v<a, ~a,b>0$  &  $\frac{a}{b}v^{\frac{1}{b} - 1}$  &  $ av^{1/b} $  &  $ ab^\eta $ \\
			Half-logistic  &  $ k\geq 0 $  &  $\frac{2k}{1-v^2}$  &  $ k\ln(\frac{1+v}{1-v}) $  &  $ 2^{-\eta}k\zeta(\eta+1) $\\
			Fr\'echet & a,b > 0 & $\frac{b^{1/a}}{av}\left(-\ln p\right)^{-\left(1+\frac{1}{a}\right)}$ & $\left(-\frac{b}{\ln v}\right)^{1/a}$& $\frac{b^{1/a}}{a\Gamma(\eta + 1)}\mathlarger{\Gamma}\left(\eta-\frac{1}{a}\right)$\\ 
			\bottomrule
			%\hline
		\end{tabular}
		\label{tab1}
		\begin{tablenotes}
			\footnotesize
			\item Note: $\zeta(s) = \sum_{k=1}^{\infty} \frac{1}{k^s}~ \text{for}~ s > 1 $ represents the Riemannian zeta function.
		\end{tablenotes}
	\end{table}
	Furthermore, there exist some lifetime models which cannot be described by any DF due to the absence of explicit expressions of DF for those models. But we can easily obtain closed form expressions for QF or qdf. For example, let us consider a random variable with qdf represented by 
	\begin{equation}
		q(v) = Kv^b(1-v)^{-(a+b)},	
		\label{2.4}
	\end{equation}  where $ a, b \in \mathcal{R} $ and $ K $ denotes the kurtosis of the distribution and hence may assume any value on the interval $ (-\infty,+\infty) $ \citep[cf. Theorem 3.2 of][]{ne}. Therefore, from the definition given by (\ref{2.i}), QFGCPE can be obtained as:
	\begin{equation}
		\mathcal{C}P\xi_Q^\eta = \begin{cases}
			\frac{K}{(b+2)^{\eta+1}} &;\text{for}~ a+b = 0 \\
			K\left(\zeta(\eta+1)- \sum_{m=1}^{b+1}\frac{1}{m^{\eta+1}}\right) & ;\text{for}~ a+b = 1 \\
			K\left[\frac{1}{(b+2)^{\eta+1}}-\frac{1}{(b+3)^{\eta+1}}\right]  & ;\text{for}~ a+b = -1 
		\end{cases}
	\end{equation} 
	\hspace*{0.2in} It is not hard to find that the random variables having the quantile function defined by (\ref{2.4}) have monotone or non-monotone reversed hazard quantile functions. Additionally, it contains several important class of distributions such as the exponential ($ a = 1, b = 0 $) generalized Pareto ($ a < 1, b = 0 $), rescaled beta ($ a > 1, b = 0 $), log-logistic ($ a = 2, b = \lambda - 1, ~ b, \lambda > 0 $) distributions and the lifetime distribution proposed by \cite{go} having the parameters $ b = \gamma - 1, a = -\gamma $ having QF given by $ Q(v) = \theta + \sigma((\gamma+1)v^r - \gamma v^{\gamma+1}) $ \citep[refer to][for more information]{ne}. 
	
\begin{figure}[]
		\includegraphics[width=0.45\textwidth]{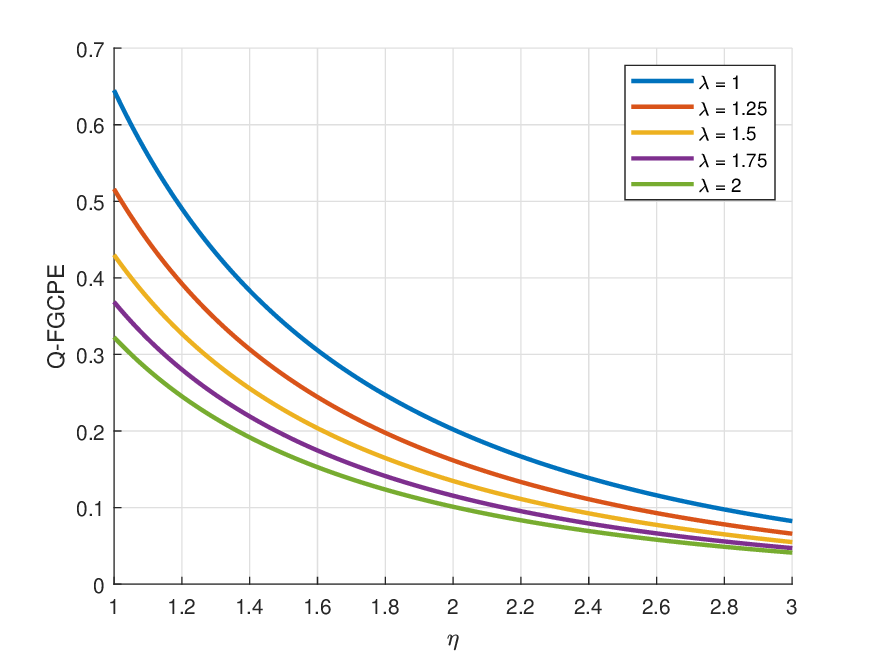} ~~~
		\includegraphics[width=0.45\textwidth]{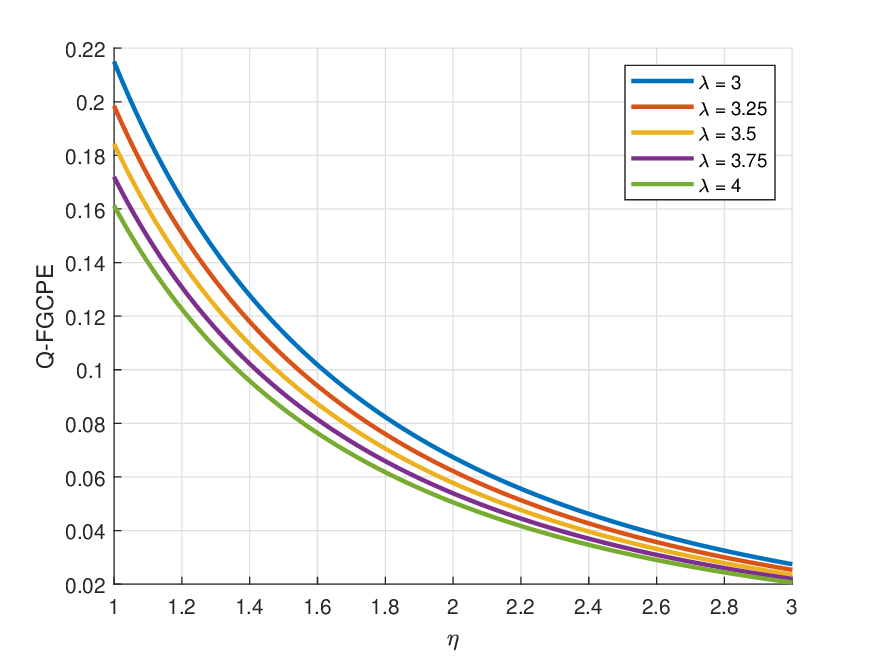}	
		\caption{Variation of QFGCPE for exponential distribution for $ \lambda = 1,1.25,1.5,1.75,2$ (left) and $ \lambda = 3,3.25,3.5,3.75,4 (right)$.} 
		\label{FIG1}
	\end{figure}

\begin{figure}[h]
		\includegraphics[width=0.45\textwidth]{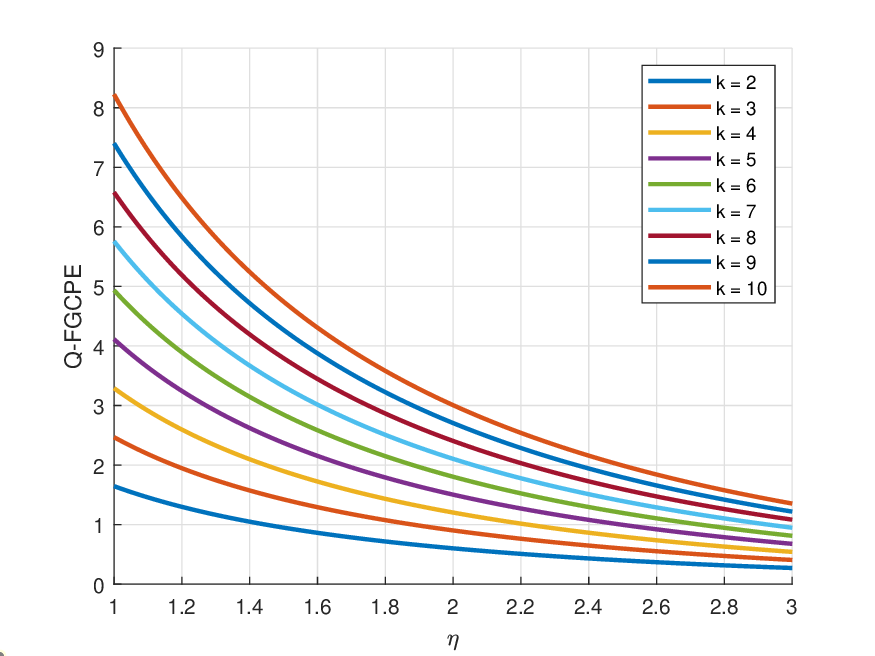} ~~~
		\includegraphics[width=0.45\textwidth]{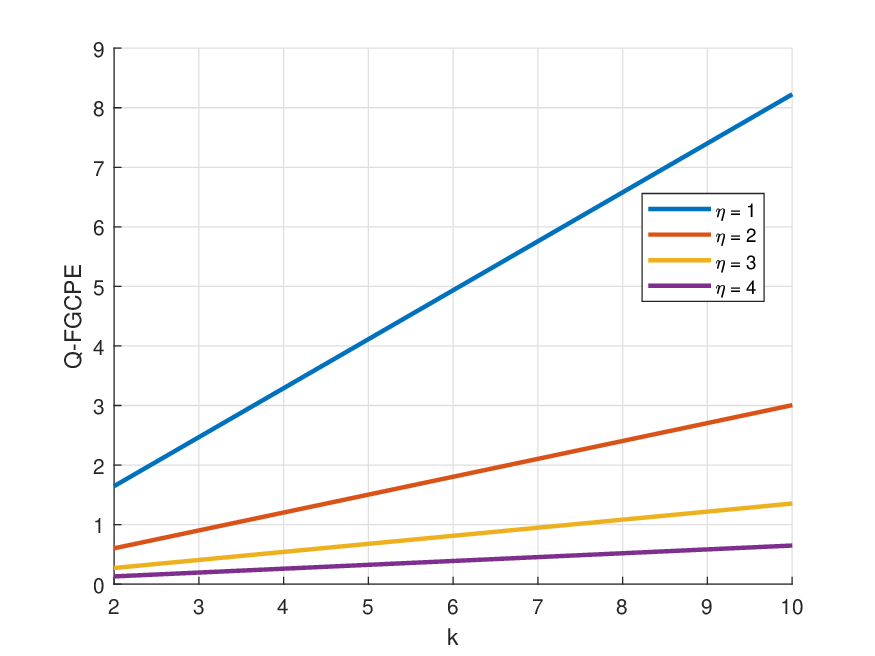}	
		\caption{Variation of QFGCPE for half-logistic distribution with respect to $ \eta $(left) and $ k $(right).} 
		\label{FIG2}
	\end{figure}

\begin{figure}[htbp]
		\begin{center}
			\includegraphics[width=0.5\textwidth]{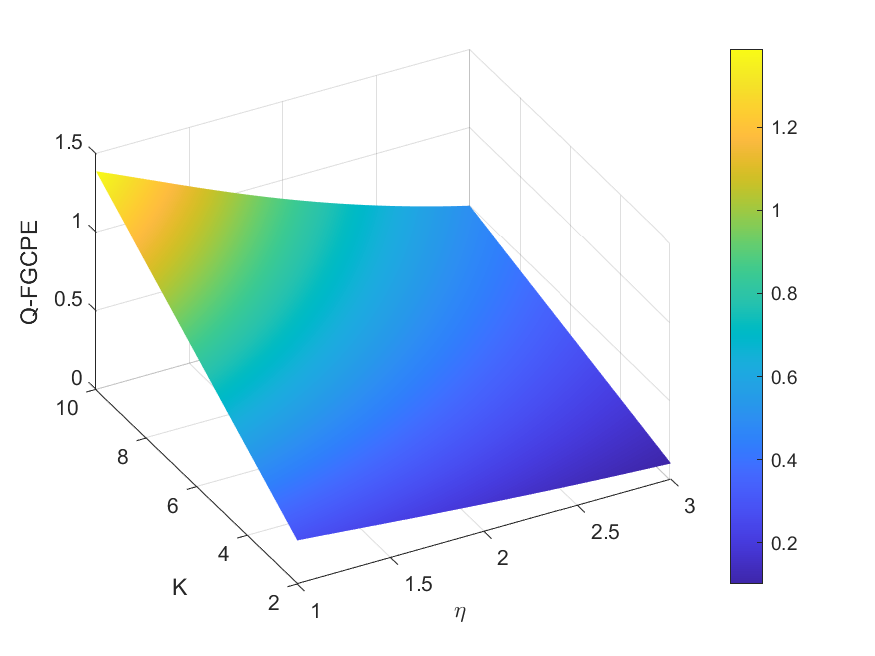}	
			\caption{Monotonic behaviour of QFGCPE for a random variable having $ q(v) = Kv^b(1-v)^{-(a+b)} $ with respect to changes in its parameters for $a=-1,~b=0$.} 
			\label{FIG3}	
		\end{center}
	\end{figure}

	Furthermore, the sensitivity of the QFGCPE function with respect to its variable parameters is portrayed through Figs. \ref{FIG1}- \ref{FIG3}. From Figs. \ref{FIG1} and \ref{FIG2}, we can observe a monotonic behavior of QFGCPE with respect to the parameters for some important classes of distributions with a tractable DF (exponential and half-logistic distribution). From Fig. \ref{FIG1}, we can observe more sensitivity of QFGCPE for lower $\lambda < 2 $ as compared to higher values $\lambda>3$ for exponential distribution. Fig. \ref{FIG2} shows that the QFGCPE function is more sensitive to the scaling parameter $k$ as compared to the fractional entropy order parameter $\eta$ for half-logistic distribution. Further, we also provide a visual representation of the changing behavior of QFGCPE with respect to the parameters $\eta$ and $K$ for the random variable with an intractable DF but having a closed form qdf given by Eq. \eqref{2.4} for the case $a+b = -1.$ through Fig. \ref{FIG3}. These illustrations highlights the importance of quantile based entropy functions in studying the variability of random variables even without an explicit DF.\\
	\hspace*{0.02in} Now, let us consider two non-negative continuous random variables $ X $ and $ Y $ with QFs $ Q_X(v) $ and $ Q_Y(v) $ respectively. Then $ Y $ is said to be the PRHM of $ X $ if and only if $ Q_Y(v) = Q_X(v^{1/\theta}) $ with qdf of $ Y ~ q_Y(v) = \frac{1}{\theta}q_X(v^{1/\theta})v^{1/\theta - 1}; ~ 0<v<1, \theta > 0 $ \citep[cf.][]{nf}. Therefore, the QFGCPE of $ Y $ is 
	\begin{equation}
		\mathcal{C}\mathcal{P}\xi_Q^\eta(Y) = \frac{1}{\theta\Gamma(\eta+1)} \int_{0}^{1} p^{\frac{1}{\theta}}(-\ln p)^\eta q_X(v^{1/\theta})dp.
		\label{2.6}
	\end{equation}  
	\begin{ex}
		Let us consider a random variables $ X $ having power distribution with qdf $ q_X(v) = \frac{a}{b}v^{\frac{1}{b}-1}$ and corresponding PRHM, $ q_Y(v) = \frac{a}{b\theta}v^{\frac{1}{b\theta}-1}, ~ a, b > 0,~ 0 < v < a < 1. $ Here, in this case, we have 
		\begin{equation*}
			\mathcal{C}\mathcal{P}\xi_Q^\eta(Y) = \frac{a(b\theta)^\eta}{(1+b\theta)^{\eta+1}}, \eta > 0.
		\end{equation*}
	\end{ex}
	From this example, one can also infer that the reversed hazard quantile functions donot show the same proportionality even though the reversed hazard rates are proportional, except for $ \theta = 1, $ where $ X $ and $ Y $ are i.i.d random variables. In other words, one can show that 
	\begin{equation*}
		R_Y(v) = b\theta v^{-1/{b\theta}} \neq \theta R_X(v) = \theta bv^{-1/b}.
	\end{equation*}
	\subsection{Some properties and bounds of QFGCPE}
	\begin{enumerate}
		\item It is straightforward to verify that Q-FGCPE is always non-negative, i.e., $ \mathcal{C}\mathcal{P}\xi_Q^\eta(Y) \geq 0. $
		\item The Q-FGCPE is dependent on the scale parameter, although it is shift-independent. Let us consider the scale and shift parameters, $ a> 0$ and $ b> 0 $, respectively. Then under affine transformation, $ Y = aX + b, $ we have $$ \mathcal{C}\mathcal{P}\xi_Q^\eta(Y) = a \mathcal{C}\mathcal{P}\xi_Q^\eta(X).$$ \\
		\textit{Proof:}
		We have from \label{1} that
		\begin{eqnarray}
			\mathcal{C}\mathcal{P}\xi^\eta(Y) &=& \frac{1}{\Gamma(\eta + 1)}\int_{0}^{\infty} G_Y(y)(-\ln G_Y(y))^\eta dy	\\
			&=& \frac{1}{\Gamma(\eta + 1)}\int_{b}^{\infty} G_X(\frac{x-b}{a})\left(-\ln G_X(\frac{x-b}{a})\right)^\eta dx, x \geq b.
			\label{2.7}
		\end{eqnarray}
		By substituting $ \frac{(x-b)}{a} $ with $ u $ in (\ref{2.7}) we get
		\begin{equation}
			\mathcal{C}\mathcal{P}\xi^\eta(Y) = \frac{a}{\Gamma(\eta + 1)}\int_{0}^{\infty} G_X(u)(-\ln G_X(u))^\eta du = a\mathcal{C}\mathcal{P}\xi^\eta(X).
			\label{2.8}
		\end{equation}
		Hence, from (\ref{2.8}), we get 
		\begin{equation*}
			\mathcal{C}\mathcal{P}\xi_Q^\eta(Y) = \frac{a}{\Gamma(\eta + 1)}\int_{0}^{1} u(-\ln u)^\eta q(u) du = a\mathcal{C}\mathcal{P}\xi_Q^\eta(X). 
		\end{equation*}
		\hfill $ \square $    
		\item Following the properties of QFs as given in \cite{nee}, if the QF corresponding to $ X, ~Q(v) = Q_1(v) + Q_2(v) $ where $ Q_1 $ and $ Q_2 $ are QFs corresponding to $ X_1 $ and $ X_2 $, respectively, then 
		\begin{equation*}
			\mathcal{C}\mathcal{P}\xi_Q^\eta(X) = \frac{1}{\Gamma(\eta + 1)}\int_{0}^{1} p(-\ln p)^\eta (q_1(p)+q_2(p))dp = \mathcal{C}\mathcal{P}\xi_{Q_1}^\eta(X) + \mathcal{C}\mathcal{P}\xi_{Q_2}^\eta(X).
		\end{equation*}
		For affirmation, we consider $ Q_1 $ to be the QF of uniform distribution and $ Q_2 $ to be the QF of uniform distribution described in Table \ref{tab1}. Then for the fixed values of the parameters: $\alpha = 0.5, ~ \lambda = 2$, we obtain 
		\[\mathcal{C}\mathcal{P}\xi_Q^\eta(X) = 1.159730,~ \mathcal{C}\mathcal{P}\xi_{Q_1}^\eta(X) = 0.353553 ~ \text{and} ~\mathcal{C}\mathcal{P}\xi_{Q_2}^\eta(X) = 0.806176. \] The results for our considered case confirm this property.
		\item If we consider two quantile functions $ Q_1 $ and $ Q_2 $ such that $ Q_1(v) + Q_2(v) = Q(v), $ then 
		\begin{equation*}
			\mathcal{C}\mathcal{P}\xi_{Q_1+Q_2}^\eta \geq \max\{\mathcal{C}\mathcal{P}\xi_{Q_1}^\eta,\mathcal{C}\mathcal{P}\xi_{Q_2}^\eta\}.
		\end{equation*}
		\textit{Proof:} This result is a straightforward consequence of properties 1 and 3.  \hfill $ \square $ 
		\item If $ Q(v) = Q_1(v)Q_2(v), $ where $ Q_1 $ and $ Q_2 $ are positive QFs, then,
		\begin{equation*}
			\mathcal{C}\mathcal{P}\xi_{Q}^\eta(X) = \frac{1}{\Gamma(\eta + 1)}\int_{0}^{1} p(-\ln p)^\eta (Q_2(p)q_1(p)+ Q_1(p)q_2(p))dp.
		\end{equation*}
		\textit{Proof:} It readily follows from the definition of QFGCPE given by (\ref{2.i}).   \hfill $ \square $ 
		\item Let $ Q_X(v) $ be the QF of $ X $. Then, the QF of $ Y = \frac{1}{X} $ is $ Q_Y(v) = \frac{1}{Q_X(1-v)}. $ Henceforth, we can express Q-FGCPE of $ Y $ as:
		\begin{equation*}
			\mathcal{C}\mathcal{P}\xi_{Q}^\eta(X) = \frac{1}{\Gamma(\eta + 1)}\int_{0}^{1} p(-\ln p)^\eta\frac{q(1-p)}{Q^2(1-p)} dp.
		\end{equation*}
		For example, let us assume that $ X $ follows Power distribution with $ Q_X(v) = v^b. $ Then, we have that $ Y = \frac{1}{X} $ has Pareto distribution with QF, $ Q_Y(v) = (1-v)^{-b}. $ Hence, the QFGCPE of $ Y $ is 
		\begin{eqnarray*}
			\mathcal{C}\mathcal{P}\xi_{Q}^\eta(Y) &=& \frac{b}{\Gamma(\eta + 1)}\int_{0}^{1} \frac{p(-\ln p)^\eta }{(1-p)^{1+b}}dp \\
			&=& \sum_{m=1}^{\infty} \frac{m}{(m+3)^{\eta + 1}} ~\text{for}~ b = 1\\
			&=& \zeta(\eta) - 3\zeta(\eta + 1) - \sum_{m=1}^{3} \left(\frac{1}{m^\eta}-\frac{3}{m^{\eta + 1}}\right) ~\text{for}~ b = 1.
		\end{eqnarray*}
		\item Let us consider a non-negative random variable $ X $ such that for $ 0 \leq \eta \leq 1, $ we have that $ \mathcal{C}\mathcal{P}\xi_{Q}^\eta(X) \leq \left[\mathcal{C}\mathcal{P}\xi_{Q}(X)\right]^\eta, $ where $ \mathcal{C}\mathcal{P}\xi_{Q}^\eta(X) $ is the quantile-based cumulative past entropy (QCPE) as defined in \cite{ssj}. This result doesnot hold when $ \eta \geq 1. $\\
		\textit{Proof:} 
		Using the relation $ p \leq p^\eta $ for $ 0 \leq p,\eta \leq 1 $ , we have %( $ p \geq p^\eta $ for $ \eta \geq 1 $)
		\begin{eqnarray*}
			\mathcal{C}\mathcal{P}\xi_{Q}^\eta(X) &=& \int_{0}^{1} p(-\ln p)^\eta q(p) dp \\
			%&\leq (\geq) & 	\int_{0}^{1} (-p\ln p)^\eta q(p) dp \\
			&\leq &  \int_{0}^{1} (-p\ln p)^\eta q(p) dp \\
			& \leq & \left[\int_{0}^{1} (-p\ln p) q(p) dp \right]^\eta ~\text{for} ~0 < \eta < 1 ~\text{(from Jensen's inequality)}\\   &=& [\mathcal{C}\mathcal{P}\xi_{Q}(X)]^\eta. 
		\end{eqnarray*} $ \hfill \square $ 
		
		We can observe that the equality condition i s satisfied when $ \eta = 1 $. 
		If we consider $g(\cdot)$ as the pdf of the random variable $ X \geq 0, $ then we obtain \[\mathcal{C}\mathcal{P}\xi_Q^\eta(X) \geq D(q)\exp [\xi_Q^\eta(X)].\]
		Here, $$ \xi_Q^\eta(X)= -\int_{0}^{1} g(Q(p))\ln g(Q(p)) dQ(p) = -\int_{0}^{1} \ln q(p) dp $$ represents the quantile-based Shannon differential entropy, and
		\[D(q) = \exp\left[\int_{0}^{1} \ln [p(-\ln p)^\eta] dp\right]\] is a non-linear function of $ \eta. $  
		\textit{Proof:} We know that $ G(x) = p = \int_{0}^{x} g(u) du. $ Then following log-sum inequality, we get \begin{eqnarray}
			\int_{0}^{1} g(Q(p))\ln \frac{g(Q(p))}{p(-\ln p)^\eta} dQ(p) &\geq& \ln \frac{1}{\int_{0}^{1} \ln [p(-\ln p)^\eta} \nonumber\\
			&=& - \ln \mathcal{C}\mathcal{P}\xi_Q^\eta(X). 
			\label{2.12}
		\end{eqnarray}  Moreover, we can rewrite the expression on the left hand side as:
		\begin{equation}
			\int_{0}^{1} g(Q(p))\ln \frac{g(Q(p))}{p(-\ln p)^\eta} dQ(p) = -\xi_Q^\eta(X) - \int_{0}^{1} \ln [p(-\ln p)^\eta] dp.
			\label{2.13}
		\end{equation} Therefore, combining the results of (\ref{2.12}) and (\ref{2.13}), we obtain
		\begin{equation}
			\ln \mathcal{C}\mathcal{P}\xi_Q^\eta(X) \geq \xi_Q^\eta(X) + \int_{0}^{1} \ln [p(-\ln p)^\eta] dp.
			\label{2.14}
		\end{equation} Thus, (\ref{2.14}) leads to the desired result given by
		\[ \mathcal{C}\mathcal{P}\xi_Q^\eta(X) \geq \exp[\xi_Q^\eta(X)] \cdot \exp\left[\int_{0}^{1} \ln [p(-\ln p)^\eta] dp\right]. \]
		\hfill $ \square $ 
	\end{enumerate}
	\begin{t1}
		Let the $QF$ and $qdf$ of a rv $X$ be denoted as $Q_X$ and $q_X$, respectively. If we consider $\varPsi > 0$ as an increasing function of $X$, then under the monotone transformation $Y = \varPsi(X), $
		\[ \mathcal{C}\mathcal{P}\xi_Q^\eta(\varPsi(X)) = \int_{0}^{1}  p(-\ln p)^\eta q_X(p)\varPsi'(Q_X(p)) dp.\] 
		\label{thm2.1}
	\end{t1}
	\textit{Proof:} Let us choose a r.v $Y$, such that $Y = \varPsi(X)$ is a positive-valued increasing function. Then we get 
	\[\mathcal{C}\mathcal{P}\xi_Q^\eta(\varPsi(X)) = \int_{0}^{1} p(-\ln p)^\eta q_Y(p)dp, ~\eta > 0.\]
	Now, we have that the pdf of Y $g_Y(y) = \frac{g_X(\varPsi^{-1}(y))}{\varPsi'(\varPsi^{-1}(y))} $ and the cdf of Y  $ G_Y(y) = G_X(\varPsi^{-1}(y)). $ This implies that $ G_Y(Q_Y(v)) = G_X(\varPsi^{-1}(Q_X(v))) \implies \varPsi^{-1}(Q_Y(v)) = G_X^{-1}(v) = Q_X(v). $ Using this information, we can express the pdf of $Y $ as:
	\[g_Y(y) = g_Y(Q_Y(v)) = \frac{g_X(\varPsi^{-1}(Q_Y(v)))}{\varPsi'(\varPsi^{-1}(Q_Y(v)))} = \frac{g_X(Q_X(v)) }{\varPsi'(Q_X(v))} = \frac{1}{q_X(v)\varPsi'(Q_X(v))}\] %= \frac{1}{q_Y(v)}.\] 
	This will lead to the final expression for the qdf of $Y $ which is the inverse of the pdf $ g_Y(y) $ given as follows:
	\[q_Y(v) = q_X(v)\varPsi'(Q_X(v)).\]
	\hfill $ \square $ 
	\begin{ex}
		Let $X$ be a r.v uniformly distributed over the interval $[0,1] $ with QF, $Q_X(v) = v $ and $q_X(v) = 1. $ Using the transformation $Y = \varPsi(X)  = X^\beta, \beta > 0 $, we get that $Y$ follows a power distribution with $Q_{\varPsi(X)(v)} = v^{1/\beta},$ such that $\varPsi(Q_X(v)) = v^\beta $ and $\varPsi'(Q_X(v)) = \beta v^{\beta-1}.$ Then the Q-FGCPE can be obtained as: \[\mathcal{C}\mathcal{P}\xi_Q^\eta(\varPsi(X)) = \frac{\beta}{(\beta + 1)^{\eta + 1}}.\]
	\end{ex}
	
	\subsection{Orderings of QFGCPE}
	The stochastic orders of quantile-based entropies can be used to compare the uncertainties of two random variables, thereby helping in management of risk \citep[cf.][and references therein]{wa}. Usually, this is done by developing a partial order relation based on the DFs of r.vs. Some fundamental stochastic order relations have been defined by \cite{sh} given as: 
	\begin{d1}
		A r.v $X$ is said to be smaller than another r.v $Y$ in the
		\begin{itemize}
			%\item[(i)] \textbf{usual stochastic order}, denoted by $X \leq_{st} Y $, if $\overline{G}_X(t) \overline{G}_Y(t) ~\forall~ t. $ 
			\item[(ii)] \textbf{hazard rate order}, denoted by $ X \leq_{hr} Y $, if $h_X(t) \geq h_Y(t) ~\forall~ t, $ where $h_X(t) = \frac{g_X(t)}{\overline{G}_X(t)}$ is the hazard rate of a r.v $X$.
			\item[(iii)] \textbf{reversed hazard rate order}, denoted by $ X \leq_{rhr} Y $, if $r_X(t) \leq r_Y(t) ~\forall~ t. $
			\item[(iv)] \textbf{dispersive order}, denoted by $X  \leq_{disp} Y $, if $G_Y^{-1}(G_X(t))-t$ in $t \geq 0.$
		\end{itemize}
		\label{def}
	\end{d1}
	Following the ordering definitions given by \cite{sh}, stated in  Definition \ref{def}, the quantile counterparts are defined as:
	\begin{d1}
		A r.v $X$ is said to be smaller than another r.v $Y$ in the
		\begin{itemize}
			\item[(v)] \textbf{hazard quantile function order}, denoted by $ X \leq_{HQ} Y $, if $H_X(v) \geq H_Y(v) ~\forall~ v \in (0,1), $ where $ H_X(v) = h(Q_X(v)) =  \frac{g(Q_X(v))}{1-v} = [(1-v)q(v)]^{-1}$ is the hazard quantile function of a r.v $X$.
			\item[(vi)] \textbf{reversed hazard rate order}, denoted by $ X \leq_{RHQ} Y $, if $R_X(v) \leq R_Y(v) ~\forall~ v \in (0,1). $ 
			\item[(iv)] \textbf{dispersive order}, denoted by $X  \leq_{disp} Y $, if  $Q_Y(v)-Q_X(v) \geq 0 ~\forall ~ v \in (0,1). $ 
		\end{itemize}	
		\label{def2.2}
	\end{d1}
	From the above defined orders, we introduce a new ordering to compare r.vs in terms of Q-FGCPE given as:
	\begin{d1}
		A r.v $X$ is said to be smaller than another r.v $Y$ in QFGCPE order, denoted as $X \leq_{QFGCPE} Y,$ if $\mathcal{C}\mathcal{P}\xi_Q^\eta(X) \leq \mathcal{C}\mathcal{P}\xi_Q^\eta(Y)$ \label{def2.3}
	\end{d1}
	Based on the QFGCPE ordering defined in Definition \ref{def2.3} and the hazard quantile function ordering given by Definition \ref{def2.2}, we can now define a property which will be followed by any two r.vs $X$ and $Y$.
	\begin{t1}
		If $X \leq_{HQ} Y $ or $X \geq_{RHQ} Y,$ then $X \leq_{QFGCPE} Y.$
	\end{t1}
	\textit{Proof:} 
	Since $X \leq_{HQ} Y \iff X \geq_{RHQ} Y $ \citep{k}, therefore, from eq. \eqref{RHQ}, we get:
	\begin{eqnarray*}
		X \leq_{HQ} Y &\implies& X \geq_{RHQ} Y \\
		&\implies& vq_X(v) \leq vq_Y(v)  \\
		&\implies& v(-\ln v)^\eta q_X(v) \leq v(-\ln v)^\eta q_Y(v) \\
		&\implies& \int_{0}^{1} v(-\ln v)^\eta q_X(v) \leq \int_{0}^{1} v(-\ln v)^\eta q_Y(v) \\
		&\implies& \mathcal{C}\mathcal{P}\xi_Q^\eta(X) \leq \mathcal{C}\mathcal{P}\xi_Q^\eta(Y).
	\end{eqnarray*}
	\hfill $ \square $ 
	
	In the following theorem, we establish a relation between all the defined orderings of Definition \ref{def2.2} and Definition \ref{def2.3} .
	\begin{t1}
		If $ X \leq_{RHQ} Y $ and $\varPsi$ is an increasing and concave function, then $X \leq_{disp} Y $ implies that $\varPsi(X) \geq_{QFGCPE} \varPsi(Y). $
		
		Similarly, if $ X \leq_{HQ} Y $ and $\varPsi$ is an increasing and convex function, then $X \leq_{disp} Y $ implies that $\varPsi(X) \leq_{QFGCPE} \varPsi(Y). $
	\end{t1}
	\textit{Proof:} Let $X \leq_{disp} Y $. This implies that $Q_X(v) \leq Q_Y(v)$. Therefore we obtain the following result for the concave and increasing function $\varPsi(X)$.
	\begin{equation}
		\varPsi'(Q_X(v)) \geq \varPsi'(Q_Y(v)) \implies 0 \leq \frac{1}{\varPsi'(Q_X(v))} \leq \frac{1}{\varPsi'(Q_Y(v))}, \label{p1}
	\end{equation}
	Thus, 
	\begin{equation}
		X\leq_{RHQ} Y \implies \frac{1}{vq_X(v)} \leq \frac{1}{vq_Y(v)}. \label{p2}
	\end{equation}
	From Eqs. \ref{p1} and \ref{p2}, we get
	\begin{equation*}
	\frac{1}{vq_X(v)\varPsi'(Q_X(v))} \leq \frac{1}{vq_X(v)\varPsi'(Q_Y(v))}	
	\end{equation*}	
	This implies
	\begin{equation}
		\int_{0}^{1} p(-\ln p)^\eta q_X(v)\varPsi'(Q_X(v)) dp \geq \int_{0}^{1} p(-\ln p)^\eta q_Y(v)\varPsi'(Q_Y(v)) dp \label{p3}
	\end{equation}
	From Eq. \eqref{p3}, we obtain
	\begin{equation*}
		\mathcal{C}\mathcal{P}\xi_Q^\eta(\varPsi(X)) \geq \mathcal{C}\mathcal{P}\xi_Q^\eta(\varPsi(Y)).
	\end{equation*}
		Similarly, the second case can be proved by following the above steps.
	\hfill $ \square $ 
	
	\section{Dynamic version of Quantile-based FGCPE}
	\begin{d1}
		The time-dependent dynamic counterpart of quantile-based fractional generalized cumulative past entropy (DQFGCPE) function can be defined as:
		\begin{equation}
			\mathcal{C}\mathcal{P}\xi_Q^\eta(X,v) = \frac{1}{v \Gamma(\eta+1)}\int_{0}^{v} p\left[\ln v - \ln p \right]^\eta q(p) dp, \eta > 0. 
			\label{3.26}
		\end{equation}
		\label{def3.1}
	\end{d1}
	When $ \eta $ is restricted to the domain $ [0,1], ~\mathcal{C}\mathcal{P}\xi_Q^\eta(X,v) $ gives the range of fractional information about the conditional probability of failure of an outcome of $ X $ upto 100u\% point of its distribution.
	
	The closed form expressions of DQFGCPE computed for some important families of distributions are provided in Table \ref{tab2}.
	
	\begin{table}[!ht]
		\centering
		\caption{Quantile-based DFGCPE for selected lifetime distributions.}
		%\begin{threeparttable}
		\renewcommand{\arraystretch}{1.3}
		\begin{tabular}{c c c c}
			\toprule
			\textbf{Distribution} & \textbf{Parameters} & $\mathbf{F(x)}$ & $\mathbf{\mathcal{C}\mathcal{P}\xi_Q^\eta(X,v)}$ \\
			\midrule
			Exponential & $\lambda > 0$ & $1-e^{-\lambda x},\; x \geq 0$ & 
			$\Gamma(\eta + 1)\varPhi(v,\eta+1,2) = \tfrac{1}{\lambda v}\big(Li_{\eta+1}(v)-v\big)$ \\[6pt]
			
			Power & $0<v<1,\; b>0$ & $x^b,\; 0 \leq x \leq l$ & 
			$vb^\eta/(b+1)^\eta $\\[6pt] %\tfrac{vb^\eta}{(b+1)^\eta}$ 
			
			Fr\'echet & $a,b > 0$ & $e^{-bx^{-a}},\; x>0$ & 
			$\tfrac{b^{1/a}}{a}(-\ln v)^{\eta-\tfrac{1}{a}}
			\,U\!\big(\eta+1,\eta+1-\tfrac{1}{a},-\ln v\big)$ \\[6pt]
			
		%	Rescaled beta & --- & --- & --- \\		
			\bottomrule
		\end{tabular}
		
		\begin{tablenotes}[flushleft]
			\footnotesize
			\item Note: \begin{enumerate}
				\item $ \varPhi(z,s,a) =  \sum_{k=0}^{\infty} \frac{z^k}{(a+k)^s} $ is the Lerch transcendent function, %\frac{1}{\Gamma(s)}\int_{0}^{\infty} \frac{t^{s-1}e^{-at}}{1-ze^{-t}} dt, $ 
				which gives the polylogarithm function for $ a = 1, $ expressed as $ Li_s(v) =  \sum_{k=1}^{\infty} \frac{v^k}{k^s}; $ 
				\item $\operatorname{U}(A,B,z) = \;=\; \frac{1}{\Gamma(A)} 
				\int_{0}^{\infty} e^{-zt} \, t^{\,A-1} (1+t)^{\,B-A-1} \, dt,$ represents the Tricomi confluent hypergeometric function, which can also be expressed in terms of Kummer’s confluent hypergeometric function $ {}_1F_1 $ as \[ U(a,b,z) \;=\; \frac{\pi}{\sin(\pi b)} \left( \frac{{}_1F_1(a,b,z)}{\Gamma(1+a-b)\Gamma(b)} - 
				\frac{z^{\,1-b}\, {}_1F_1(1+a-b,2-b,z)}{\Gamma(a)\Gamma(2-b)} \right). \]
			\end{enumerate} 
		\end{tablenotes}
		%	\end{threeparttable}
	\label{tab2}
\end{table}
\begin{figure}[htbp]
	\begin{center}
		\includegraphics[width=0.7\textwidth]{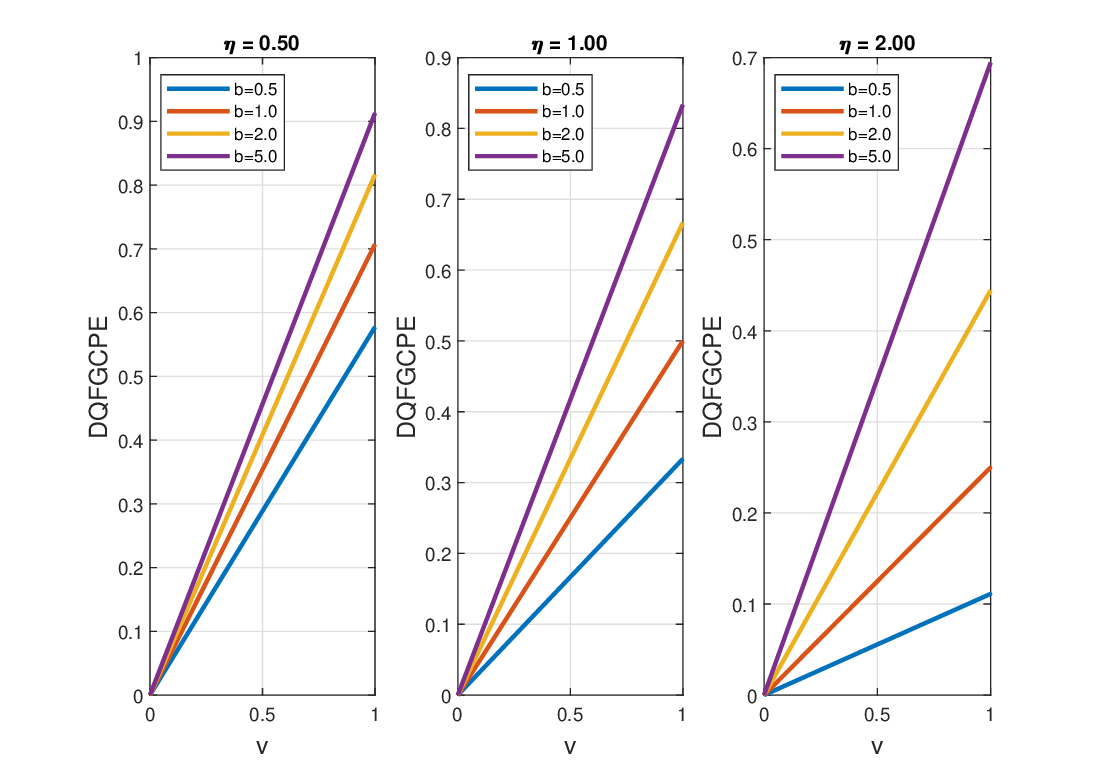}	
		\caption{Monotonic behaviour of DQFGCPE for power distribution with respect to $v$.} 
		\label{FIG4}	
	\end{center}
\end{figure} 
\newpage
Fig. \ref{FIG4} depicts the increase in sensitivity of the DQFGCPE measure with increase in the entropy order parameter $\eta$ values for the power distribution. We also see that the dynamic FGCPE follow an increasing trend with respect to $v$ for all values of the distribution parameter $b > 0$ and the entropy parameter $\eta > 0$.
The following theorem discusses about the DQFGCPE expression under monotone transformation $Y = \varPsi(X)$, with $\varPsi$ being an increasing function.
\begin{t1}
	Let the $QF$ and $qdf$ of a rv $X$ be denoted as $Q_X$ and $q_X$, respectively. If we consider $\varPsi > 0$ as an increasing function, then DQFGCPE can be expressed as:
	\[  \mathcal{C}\mathcal{P}\xi_Q^\eta(\varPsi(X),v) = \frac{1}{v\Gamma(\eta + 1)}\int_{0}^{1} p(\ln u - \ln p)^\eta q_X(p)\varPsi'(Q_X(p)) dp. \] 
\end{t1}
\textit{Proof:} This result directly follows from the property of $\mathcal{C}\mathcal{P}\xi_Q^\eta(\varPsi(X))$ proved in Theorem \ref{thm2.1}.
\hfill $ \square $ 

On the basis of DQFGCPE, two non-parametric classes of life distributions are defined as follows:
\begin{d1}
	A r.v $X$ is said to have increasing (or decreasing) DQFGCPE, denoted as IDQFGCPE (or DDQFGCPE) only if $\mathcal{C}\mathcal{P}\xi_Q^\eta(X,v)$ is increasing (or decreasing) in $v \geq 0 $.
\end{d1}
In other words, $\frac{d}{dv}\mathcal{C}\mathcal{P}\xi_Q^\eta(X,v) \geq (or \leq) 0 \implies \frac{d}{d\eta}\mathcal{C}P\xi_Q^\eta(X,v) \geq (or \leq) $ whenever $X$ is IDQFGCPE (or DDQFGCPE).

\begin{d1}
	If $\mathcal{C}P\xi_Q^\eta(X,v) \leq \mathcal{C}\mathcal{P}\xi_Q^\eta(Y,v)$, then the r.v $X$ is said to have lesser DQFGCPE than the r.v $Y$, represented as $X \leq_{QDFGCPE} Y $. 
\end{d1}

For example, let $X$ and $Y$ be exponentially distributed r.vs with mean failure rates $\lambda_1$ and $\lambda_2,$ respectively, then $\lambda_1 \geq \lambda_1 $ implies that $ X \leq_{QDFGCPE} Y $.

Based on the above two new classes of life distributions in terms of DQFGCPE, we discuss some monotonic and stochastic ordering properties of DQFGCPE in the following theorems.
\begin{t1}
	Let us consider a r.v $Y = \varPsi(X)$ transformed from the non-negative r.v $X,$ where $\varPsi(\cdot)$ is an increasing real-valued, positive and convex (or concave) function. Then, for $\eta > 0,$ $Y$ is IDQFGCPE (or DDQFGCPE) if $X$ is IDQFGCPE (or DDQFGCPE).
	\label{thm3.2}
\end{t1}
%\begin{proof}
\textit{Proof:} According to the definition of DQFGCPE given by Eq. \eqref{3.26}, we obtain
	\begin{eqnarray*}
	  \mathcal{C}P\xi_Q^\eta(Y,v) &=& \frac{1}{v \Gamma(\eta+1)}\int_{0}^{v} p\left[\ln v - \ln p \right]^\eta q_Y(p) dp, \qquad \eta > 0 \\
	  &=& \frac{1}{v \Gamma(\eta+1)}\int_{0}^{v} p\left[\ln v - \ln p \right]^\eta q_X(p)\varPsi'(Q_X(p)) dp, \eta > 0. 
	\end{eqnarray*}
Considering the increasing, non-negative and convex (or concave) nature of $\varPsi$, we can claim $\varPsi'(Q_X(v))$ to follow a similar increasing (or decreasing) and non-negative behavior with respect to increase in $v$. Hence following Theorem 3.2 of \cite{ses}, we get that $ \mathcal{C}P\xi_Q^\eta(Y,v)$ is increasing (or decreasing) in $v$. Hence, for $\eta > 0,~ \mathcal{C}P\xi_Q^\eta(X,v) $ is increasing (or decreasing) in $v$ implies that $\mathcal{C}P\xi_Q^\eta(Y,v) $ is also increasing (or decreasing) in $v$. \hfill $ \square $ 

This property can be justified by taking the case of affine transformation $Y = aX + b, ~ a > 0, ~ b \geq 0.$
%\end{proof}

\begin{t1}
	If $X \leq_{disp} Y, $ then $X \leq_{DQFGCPE} Y.$
\end{t1}
\textit{Proof:} Following the definition of dispersive order from Definition \ref{def2.2}, we get that $q_X(v) \leq q_Y(v).$ Hence, for $X \leq_{disp} Y, $ we obtain
\begin{eqnarray*}
	\mathcal{C}P\xi_Q^\eta(X,v) &=& \frac{1}{v \Gamma(\eta+1)}\int_{0}^{v} p\left[\ln v - \ln p \right]^\eta q_X(p) dp \\
	&\leq& \frac{1}{v \Gamma(\eta+1)}\int_{0}^{v} p\left[\ln v - \ln p \right]^\eta q_Y(p) dp~~ = ~~\mathcal{C}P\xi_Q^\eta(Y,v).  
\end{eqnarray*} \hfill $ \square $

\begin{t1}
	Let us consider that two r.vs $X$ and $Y$ following the ordering $X \leq_{DQFGCPE} Y.$ Then this ordering property is preserved for an increasing non-negative convex function $\varPsi.$ In other words, we have $X \leq_{DQFGCPE} Y~\implies ~ \varPsi(X) \leq_{DQFGCPE} \varPsi(Y).$
	\label{thm3.4}
\end{t1}
\textit{Proof:} The proof follows from %Definition \ref{def3.1} and 
the results of Theorems \ref{thm3.2} and \ref{thm3.4} and the property that $(\ln v - \ln p)> 0 $ if $p < v.$
%\newpage
\section{Non-parametric estimation of QFGCPE}
In this section, we develop a non-parametric estimator for the quantile-based fractional generalized cumulative past entropy (QFGCPE). The proposed estimator is constructed directly from the sample quantiles, thereby avoiding assumptions about the underlying distributional form. Let $X_1, X_2, \dots, X_n$ be a sequence of $n$ independent and identically distributed (i.i.d.) random variables with cumulative distribution function (cdf) $F(x)$ and corresponding quantile function (QF) $Q(v)$. We denote their order statistics by  
\[
X_{1:n} \leq X_{2:n} \leq \cdots \leq X_{n:n},
\]  
where $X_{k:n}$ represents the $k$-th smallest observation in the sample, for $k = 1, 2, \dots, n$.

The order statistics naturally give rise to the empirical distribution function (edf), defined as  
\[
F_n(x) = \frac{1}{n} \sum_{i=1}^n \mathbf{1}_{\{X_i \leq x\}}, \qquad x \in \mathbb{R},
\]  
which serves as a non-parametric estimator of the true distribution function $F(x)$. Correspondingly, the empirical quantile function (EQF) is obtained by inverting $F_n(x)$, and is formally given by  
\[
Q_n(v) = \inf \{x \in \mathbb{R} : F_n(x) \geq v \}, \qquad 0 < v < 1.
\]  

In particular, for $v \in \big(\tfrac{k-1}{n}, \tfrac{k}{n}\big]$, the EQF $Q_n(v)$ coincides with the $k$-th order statistic $X_{k:n}$ \citep{p}. A smoothed form of the estimator of EQF is defined as \citep{ses}:
\begin{equation*}
	\overline{Q}_n(v) = n\left(\tfrac{k}{n}-v\right) X_{k-1:n} +  n\left(v-\tfrac{k-1}{n}\right) X_{k:n}, \qquad v \in \big(\tfrac{k-1}{n}, \tfrac{k}{n}\big),~k = 1,2,\dots,n.
\end{equation*}

Consquently, the empirical quantile density function (eqdf) can be obtained as:
\begin{equation}
	\overline{q}_n(v) = n(X_{k:n}-X_{k-1:n}). % \qquad v \in \big(\tfrac{k-1}{n}, \tfrac{k}{n}\big).
	\label{4.1}
\end{equation}

Therefore, from Eqs. \eqref{4.1} and \eqref{2.i}, we compute the non-parametric estimator of Q-FGCPE as:
\begin{equation}
	\widehat{\mathrm{\mathcal{C}\mathcal{P}\xi}}_Q^\eta(X) =  \frac{1}{\Gamma(\eta+1)}\int_{0}^{1} p(-\ln p)^\eta \overline{q}_n(v) dp,  \eta > 0
	\label{4.2}
\end{equation}  
The final expression, derived from approximating the integral in Eq. \eqref{4.2} by summation over n intervals of $v \in (0,1)$, can be written as:
\begin{equation}
	\widehat{\mathrm{\mathcal{C}\mathcal{P}\xi}}_Q^\eta(X) =  \frac{1}{\Gamma(\eta+1)}\sum_{k=1}^{n} F_n(X_{k:n})[-\ln F_n(X_{k:n})]^\eta ~n(X_{k:n}-X_{k-1:n}) (S_{k:n}-S_{k-1:n}),
	\label{4.3}
\end{equation}
where $F_n(X_{k:n})$ is the edf and $S_{k:n}$ is defined as:
\[ S_{k:n}
\begin{cases}
	0, \quad k = 0 \\
	F(X_{k:n}) = \tfrac{k}{n}, \quad k = 1,2,\dots,n-1 \\
	1, \quad k = n 
\end{cases} \]

Hence, Eq. \eqref{4.3} can be rewritten as:
\begin{equation}
	\widehat{\mathrm{\mathcal{C}\mathcal{P}\xi}}_Q^\eta(X) =  \frac{1}{\Gamma(\eta+1)}\sum_{k=1}^{n-1} \frac{k}{n}\left(-\ln \frac{k}{n}\right)^\eta (X_{k:n}-X_{k-1:n}).
	\label{4.4}
\end{equation}  
\section{Simulation Study}

To evaluate the finite-sample performance of the proposed QFGCPE estimator defined in Eq. \eqref{4.4} under the standard exponential and Govindarajalu distributions, we conducted a Monte Carlo simulation experiment. Independent random samples of sizes $n = 50, 100, 500, 1000,$ and $5000$ are generated, with each setting replicated $N_{\text{sim}} = 500$ times. For each dataset, the QFGCPE is estimated, and its empirical properties (bias and MSE) are compared with the corresponding theoretical values computed from \ref{2.i} (see Tables \ref{simtab1} - \ref{simtab4}). 

In our simulation study, the Govindarajalu distribution (Govindarajalu, 1977) is deliberately selected to highlight the role of the quantile-based definition of FGCPE, as it accommodates random variables with intractable distribution functions \citep{ses}. Furthermore, this quantile entropy has proven useful in modeling bathtub-shaped lifetime data \citep{nsb}. The QF of the Govindarajalu distribution is given by:
\begin{equation}
	Q(v) = \alpha + \beta((\gamma+1)v^\gamma - \gamma v^{\gamma+1}), \alpha \in (-\infty,+\infty), \beta,\gamma > 0,~ v \in (0,1).
\end{equation}

\begin{table}[]%[htbp]
	\centering
	\caption{Comparison of Empirical vs Theoretical QFGCPE ($\eta=0.5$) for Exponential($\lambda = 1$) Distribution}
	\renewcommand{\arraystretch}{1.3}
	\begin{tabular}{c c c c c c}
		\toprule
		$\mathbf{n}$ & \textbf{Mean Empirical} & \textbf{Bias} & \textbf{MSE} & \textbf{RMSE} & \textbf{Theoretical Value} \\
		\midrule
		50   & 1.3861  & $-0.2262$ & 0.1364   & 0.3694 & 1.6124 \\
		100  & 1.4610  & $-0.1514$ & 0.0761  & 0.2759 & 1.6124 \\
		500  & 1.5326  & $-0.0797$ & 0.0222 & 0.1490 & 1.6124 \\
		1000 & 1.5573  & $-0.0550$ & 0.0120  & 0.1094 & 1.6124 \\
		5000 & 1.5883  & $-0.0241$ & 0.0025 & 0.0499 & 1.6124 \\
		\bottomrule
	\end{tabular}
	\label{simtab1}
\end{table}

\begin{table}[]%[htbp]
	\centering
	\caption{Comparison of Empirical vs Theoretical QFGCPE ($\eta=0.75$) for Exponential($\lambda = 1$) Distribution}
	\renewcommand{\arraystretch}{1.3}
	\begin{tabular}{c c c c c c}
		\toprule
		$\mathbf{n}$ & \textbf{Mean Empirical} & \textbf{Bias} & \textbf{MSE} & \textbf{RMSE} & \textbf{Theoretical Value} \\
		\midrule
		50   & 0.91716  & $-0.04516$ & 0.02760   & 0.1661 & 0.96232 \\
		100  & 0.94243  & $-0.01989$ & 0.01356   & 0.1164 & 0.96232 \\
		500  & 0.95234  & $-0.00998$ & 0.00330   & 0.0574 & 0.96232 \\
		1000 & 0.95623  & $-0.00609$ & 0.00164   & 0.0406 & 0.96232 \\
		5000 & 0.96123  & $-0.00109$ & 0.00028   & 0.0169 & 0.96232 \\
		\bottomrule
	\end{tabular}
	\label{simtab2}
\end{table}

\begin{table}[]%[ht]
	\centering
	\caption{Comparison of Empirical vs Theoretical QFGCPE for Govindarajalu($\alpha,\beta,\gamma$) distribution with $(\alpha,\beta,\gamma) = (1,2,2)$ and $\eta=0.25$}
	\renewcommand{\arraystretch}{1.3}
	\begin{tabular}{c c c c c c}
		\toprule
		$\mathbf{n}$ & \textbf{Mean Empirical} & \textbf{Bias} & \textbf{MSE} & \textbf{RMSE} & \textbf{Theoretical Value} \\
		\midrule
		50   & 0.9045  & $-1.35 \times 10^{-2}$ & $3.80 \times 10^{-3}$ & $6.16 \times 10^{-2}$ & 0.91802 \\
		100  & 0.9148  & $-3.20 \times 10^{-3}$ & $1.84 \times 10^{-3}$ & $4.29 \times 10^{-2}$ & 0.91802 \\
		500  & 0.9171  & $-9.40 \times 10^{-4}$ & $3.78 \times 10^{-4}$ & $1.94 \times 10^{-2}$ & 0.91802 \\
		1000 & 0.9176  & $-4.20 \times 10^{-4}$ & $1.84 \times 10^{-4}$ & $1.36 \times 10^{-2}$ & 0.91802 \\
		5000 & 0.9182  & $ 1.53 \times 10^{-4}$ & $3.55 \times 10^{-5}$ & $5.96 \times 10^{-3}$ & 0.91802 \\ \bottomrule
	\end{tabular}
	\label{simtab3}
\end{table}

\begin{table}[]%[ht]
	\centering
	\caption{Comparison of Empirical vs Theoretical QFGCPE for Govindarajalu($\theta,\sigma,\beta$) distribution with $(\alpha,\beta,\gamma) = (1,2,2)$ and $\eta=0.75$}
	\renewcommand{\arraystretch}{1.3}
	\begin{tabular}{c c c c c c}
		\toprule
		$\mathbf{n}$ & \textbf{Mean Empirical} & \textbf{Bias} & \textbf{MSE} & \textbf{RMSE} & \textbf{Theoretical Value} \\
		\midrule
		50   & 0.68147 & $-1.26 \times 10^{-2}$ & $9.35 \times 10^{-4}$ & $3.06 \times 10^{-2}$ & 0.69411 \\
		100  & 0.68875 & $-5.35 \times 10^{-3}$ & $3.87 \times 10^{-4}$ & $1.97 \times 10^{-2}$ & 0.69411 \\
		500  & 0.69315 & $-9.53 \times 10^{-4}$ & $7.51 \times 10^{-5}$ & $8.67 \times 10^{-3}$ & 0.69411 \\
		1000 & 0.69320 & $-9.06 \times 10^{-4}$ & $3.83 \times 10^{-5}$ & $6.19 \times 10^{-3}$ & 0.69411 \\
		5000 & 0.69407 & $-3.03 \times 10^{-5}$ & $6.45 \times 10^{-6}$ & $2.54 \times 10^{-3}$ & 0.69411 \\ \bottomrule
	\end{tabular}
	\label{simtab4}
\end{table}

Tables~\ref{simtab1}-\ref{simtab4} summarize the results in terms of the empirical mean of the estimator, bias, mean squared error (MSE), root mean squared error (RMSE), and the theoretical benchmark. As expected, the estimator obeys the asymptotic property and exhibits clear consistency across all the chosen distributions: both the bias and MSE decrease systematically with increasing sample size $n$, and the empirical mean converges towards the theoretical value. For small samples ($n=50$), the bias and RMSE are relatively larger, but these values diminish rapidly as $n$ grows. At $n=5000$, the bias is essentially negligible and the RMSE is low, indicating high accuracy.

In addition to point estimation, we also assessed the interval estimation performance using bootstrap-based confidence intervals. For each simulated dataset, $B_{\text{boot}} = 500$ bootstrap replications are drawn, and 95\% confidence intervals are constructed using the percentile method. The empirical coverage probability and the Monte Carlo standard error (MCSE) are then evaluated across replications. The results for standard exponential and Govindarajalu(1,2,2) distributions with $\eta = 0.75$ are displayed in Tables \ref{simtab5} and \ref{simtab6}, respectively. Table \ref{simtab5} demonstrate that the coverage improves with larger sample size $n$ for standard exponential distribution. For $n=50$, the coverage is about 75.8\%, but it rises to over 90\% for $n=500$ and beyond, with MCSE values around 0.013--0.019. We observe similar trends with Govindarajalu(1,2,2) distribution as shown in Table \ref{simtab6}. The coverage probabilities improve with increasing sample size, approaching the nominal $0.95$ level for $n \geq 500$. For smaller samples ($n=50,100$), coverage is slightly lower, reflecting finite-sample bias. The decreasing MCSE values with larger $n$ indicate greater precision and stability of the coverage estimates.
This confirms that bootstrap confidence intervals provide reliable uncertainty quantification in moderate to large samples, though slightly conservative intervals may be needed for small samples. Overall, the results suggest that the proposed estimator yields reliable coverage properties, particularly in moderate to large samples. 

\begin{table}[h!]
	\centering
	\caption{Bootstrap coverage probabilities and Monte Carlo standard errors (MCSE) for 95\% percentile confidence intervals of QFGCPE for Exponential($\lambda=1$) distribution with $\eta = 0.75$.}
	\renewcommand{\arraystretch}{1.3}
	\begin{tabular}{c c c}
		\toprule
		$\mathbf{n}$ & \textbf{Coverage Probability} & \textbf{MCSE} \\
		\midrule
		50   & 0.758 & 0.019 \\
		100  & 0.846 & 0.016 \\
		500  & 0.912 & 0.013 \\
		1000 & 0.924 & 0.013 \\
		%5000 & 0.938 & 0.011 \\
		\bottomrule
	\end{tabular}
	\label{simtab5}
\end{table}

\begin{table}[h!]
	\centering
	\caption{Bootstrap coverage probabilities and Monte Carlo standard errors (MCSE) for 95\% percentile confidence intervals of QFGCPE for Govindarajalu(1,2,2) distribution with $\eta = 0.75$.}
	\renewcommand{\arraystretch}{1.3}
	\begin{tabular}{c c c}
		\toprule
		$\mathbf{n}$ & \textbf{Coverage Probability} & \textbf{MCSE} \\
		\midrule
		50   & 0.892 & 0.0139 \\
		100  & 0.924 & 0.0119 \\
		500  & 0.958  & 0.0090 \\
		1000 & 0.942 & 0.0105 \\
		%5000 & 0.938 & 0.011 \\
		\bottomrule
	\end{tabular}
	\label{simtab6}
\end{table}

\section{Validity and Sensitivity of QFGCPE with Logistic Maps} 	
We verify the validity and sensitivity of our proposed QFGCPE measure by analyzing its behavior for both periodic and chaotic regimes through simulation studies on a logistic map defined with the help of a canonical model of deterministic chaos given as: %recurrence relation 
\begin{equation}
	x_{n+1} = cx_n(1-x_n), \qquad x \in [0,1].
	\label{6.1}
\end{equation} This model \eqref{6.1} generates sequences that switch between stable, periodic and chaotic behaviors depending on the parameter $c.$ By varying $c$ in the range $[0,4],$ one can generate a time series of predictable (low entropy) or chaotic (high entropy) nature. For $c > 3, $ the system is expected to exhibit chaotic behavior while a stable and periodic behavior is displayed by the system for $c$ values lesser than $3$. Here, we conduct this study by choosing an initial value of $x_0 = 0.1$ and assigning different values to the control parameter $c \in [0,4]$. The selected values of $c$ for the validation of the proposed entropy estimator are $c = 1, 1.5, 2, 2.5, 3, 3.5, 4.$   

\begin{figure}[h]
	\centering
	\includegraphics[width=0.7\textwidth]{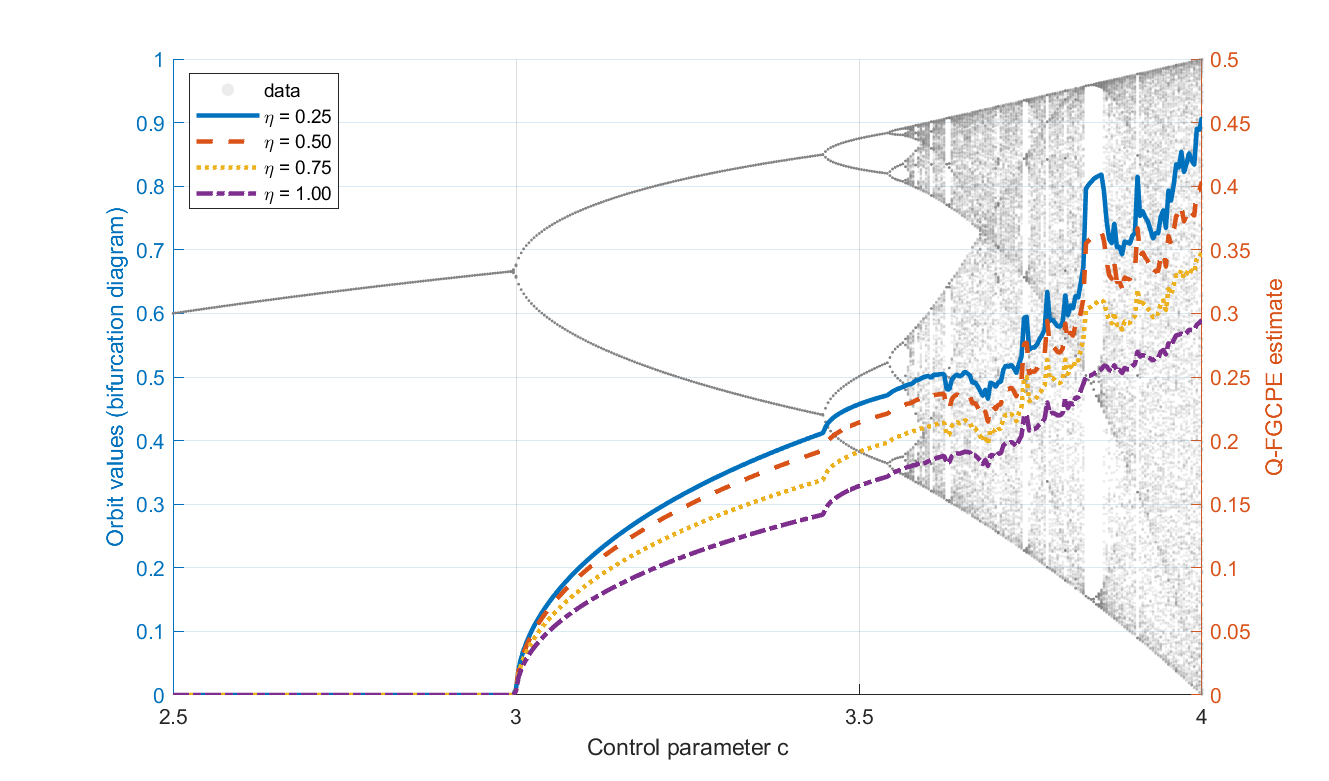}
	\caption{Logistic map: bifurcation vs QFGCPE with respect to $c$ for different $ \eta $ values.} 
	\label{Valid2}
\end{figure} 
\begin{figure}[h]
	\centering
	\includegraphics[width=0.7\textwidth]{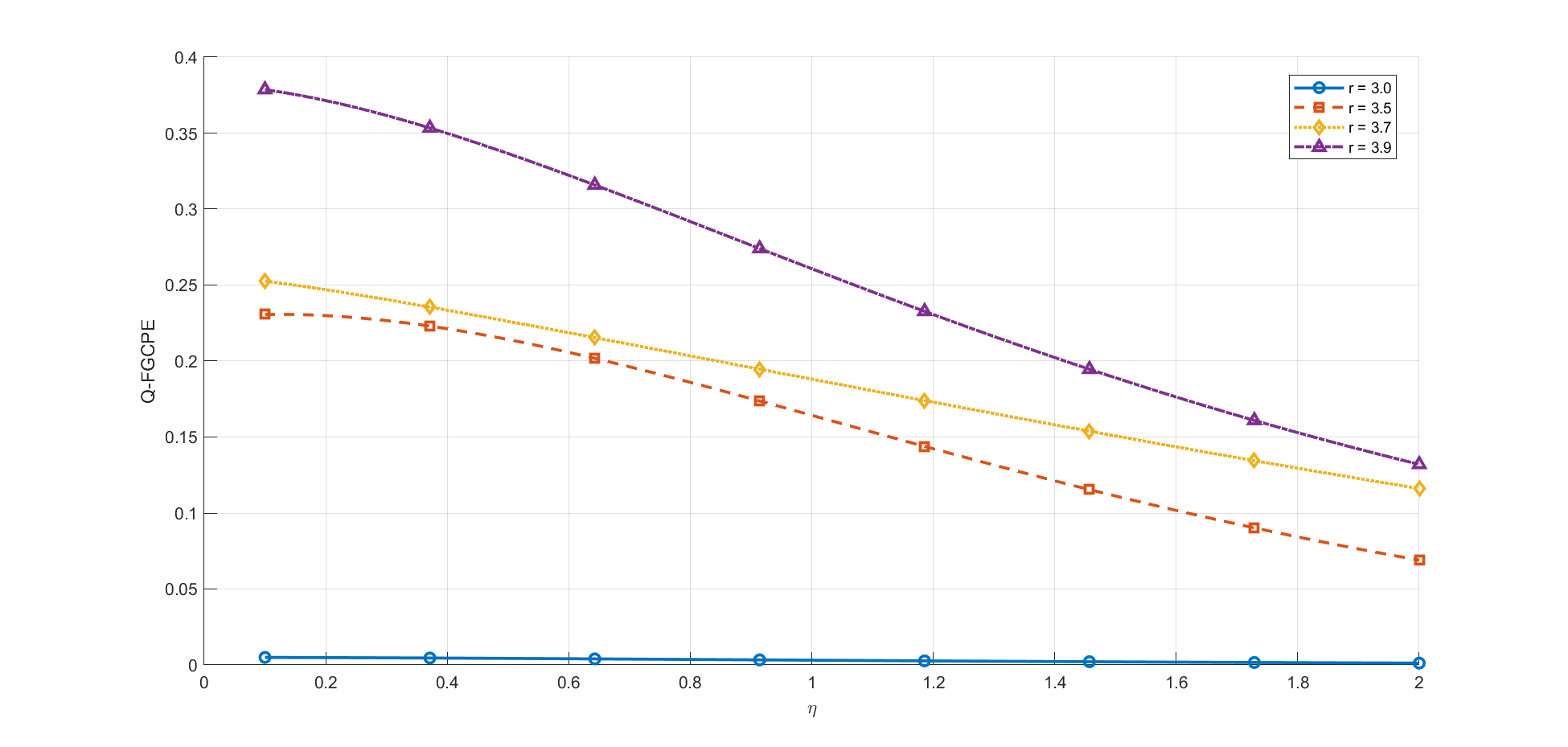} 	
	\caption{Logistic map: QFGCPE with respect to $ \eta $(left) for different $c$ values.} 
	\label{Valid1}
\end{figure} 

%\textbf{Variation across control parameter $r$.}  
When QFGCPE is plotted against the logistic map parameter $c \in [1,4]$ for different fractional orders $\eta$ as shown in Fig. \ref{Valid2}, the estimator remains close to zero in the stable regime ($c \lesssim 3$), where the dynamics converge to fixed points or low-period cycles and randomness is minimal. As $c$ increases beyond the bifurcation threshold, QFGCPE rises sharply, mirroring the onset of complex dynamics and chaos. The growth pattern is consistent with the transition structure observed in the bifurcation diagram, %and correlates positively with the Lyapunov exponent, 
thereby validating that the measure captures the degree of dynamical uncertainty. Distinct $\eta$ values yield qualitatively similar trends but different scales, with smaller $\eta$ giving more pronounced sensitivity to small fluctuations, while larger $\eta$ smooths variability.  

%\textbf{Variation across fractional order $\eta$.}  
Fixing $c$ at representative values (periodic window $c=3.2$, chaotic regime $c=3.7$), we examined QFGCPE as a function of $\eta$ illustrated by Fig. \ref{Valid1}. The estimator increases monotonically with $\eta$, but the rate and curvature differ depending on the dynamical regime. For low-period orbits represented by lower $c$ values, the sensitivity to $\eta$ is weak, consistent with low dynamical complexity. In contrast, for chaotic regimes, the growth of QFGCPE with $\eta$ is steeper, indicating stronger responsiveness of the measure to fractional order. 

Cumulatively, these evidences point towards the argument that the QFGCPE estimator is both valid, reflecting the transition from regular to chaotic dynamics, and sensitive enough to distinguish levels of complexity across both system parameter $c$ and entropy order $\eta$. This dual perspective highlights QFGCPE as a flexible and robust tool for quantifying dynamical uncertainty.  	

\section{Conclusion}
The overall contributions of this work are threefold:
\begin{itemize}
	\item[(i)] the introduction of a quantile-based analogue of fractional generalized cumulative entropy and its theoretical properties;
	\item[(ii)] explicit derivations for important lifetime distributions and extensions to dynamic settings; and
	\item[(iii)] the construction and validation of a nonparametric estimator, supported by simulation and dynamical system comparisons.
\end{itemize} 
These developments place QFGCPE within the growing family of quantile-based entropy measures, enriching the statistical toolbox for uncertainty quantification in reliability, survival analysis, and dynamical modeling.  

\section*{Data availability statement}
All data used in this study are electronically available and are all cited in the bibliography section.
\section*{Disclosure statement}
On behalf of all authors, the corresponding author states that there is no conflict of interest.		

\end{document}